\documentclass[sn-mathphys-num]{sn-jnl}

\usepackage{graphicx}%
\usepackage{multirow}%
\usepackage{amsmath,amssymb,amsfonts,amsthm,latexsym}%
\usepackage{amsthm}%
\usepackage{mathrsfs}%
\usepackage[title]{appendix}%
\usepackage{xcolor}%
\usepackage{textcomp}%
\usepackage{manyfoot}%
\usepackage{booktabs}%
\usepackage{algorithm}%
\usepackage{algorithmicx}%
\usepackage{algpseudocode}%
\usepackage{listings}%
\usepackage{textcomp}
\usepackage{tikz-cd}
\usepackage[all]{xy}
\usepackage{mathpazo}

\usepackage{comment}

\theoremstyle{thmstyleone}%
\newtheorem{theorem}{Theorem}%
\newtheorem{proposition}[theorem]{Proposition}%
\newtheorem{lemma}[theorem]{Lemma}%
\newtheorem{corollary}[theorem]{Corollary}%

\theoremstyle{thmstyletwo}%
\newtheorem{example}{Example}%
\newtheorem{remark}{Remark}%

\theoremstyle{thmstylethree}%
\newtheorem{definition}{Definition}%

\newcommand\set[1]{\mbox{$\{#1\}$}}
\newcommand\cat[1]{\mbox{$\mathbf{#1}$}}
\newcommand\ds[1]{\mbox{$\displaystyle{#1}$}}
\newcommand\map[3]{\mbox{$#1:#2\rightarrow #3$}}

\newcommand\trip[1]{\mbox{$\mathbb{#1}$}}

\newcommand{\id}{\mbox{$\mathfrak{I}$}}

\newcommand\comp[2]{\mbox{$#1\cdot #2$}}
\newcommand\ocomp[2]{\mbox{$#1\circ #2$}}
\newcommand\join[2]{\mbox{$#1\vee #2$}}
\newcommand\meet[2]{\mbox{$#1\wedge #2$}}
\newcommand\mon[3]{\mbox{$\left(#1,#2,#3\right)$}}
\newcommand\fak[1]{\mbox{$#1^{\varphi}$}}
\newcommand\powcat[2]{\mbox{$#1^{\mathbb{#2}}$}}
\newcommand\alg[1]{\mbox{$\mathbf{Alg}\left(\mathbb{#1}\right)$}}
\newcommand\algn[2]{\mbox{$\mathbf{Alg}\left(\mathbb{#1}_{#2}\right)$}}
\newcommand\coalg[1]{\mbox{$\mathbf{Coalg}\left(\mathbb{#1}\right)$}}
\newcommand\coalgn[2]{\mbox{$\mathbf{Coalg}\left(\mathbb{#1}_{#2}\right)$}}

\newcommand\Adj[4]{\mbox{$$\xymatrix{#1 \ar@<1ex>[rr]^{#2}_*-<0.55ex>{_\perp} & & \ar@<1ex>[ll]^{#3} #4}$$}}
\newcommand\eq[6]{\mbox{$$\xymatrix{#1\ar[r]^{#2} & #3\ar@<0.5ex>[r]^{#4} \ar@<-0.5ex>[r]_{#5}  & #6}$$}}

\newcommand{\arrowtcupp}[2]{\arrow[bend left=50, ""{name=U, below,inner sep=1}]{#1}\arrow[Rightarrow,from=U,to=MU,"#2"]}
\newcommand{\arrowtclow}[2]{\arrow[bend right=50, ""{name=L,inner sep=1}]{#1}\arrow[Rightarrow,from=LM,to=L]{}[]{#2}} 
\newcommand{\arrowtcmid}[2]{\arrow[""{name=MU,inner sep=1},""{name=LM,below,inner sep=1}]{#1}[pos=.1]{#2}}

\begin{document}

\title[Monadic aspects of the ideal lattice functor on the category of distributive lattices]{Monadic aspects of the ideal lattice functor on the category of distributive lattices}


\author*[1]{\fnm{Ando} \sur{Razafindrakoto}}\email{arazafindrakoto@uwc.ac.za}

\affil[1]{\orgdiv{Department of Mathematics and Applied Mathematics}, \orgname{University of the Western Cape}, \orgaddress{\street{Robert Sobukwe Rd}, \city{Bellville}, \postcode{7535}, \state{Cape Town}, \country{South Africa}}}

\abstract{It is known that the construction of the frame of ideals from a distributive lattice induces a monad whose algebras are precisely the frames and frame homomorphisms. Using the Fakir construction of an idempotent approximation of a monad, we extend B. Jacobs' results on lax idempotent monads and show that the sequence of monads and comonads generated by successive iterations of this ideal functor on its algebras and coalgebras do not strictly lead to a new category. We further extend this result and provide a new proof of the equivalence between distributive lattices and coherent frames by showing that when the first inductive step in the Fakir construction is the identity monad, then the ambient category is equivalent to the free algebras.}

\keywords{Monad, algebras, distributive lattices, frames, ideal functor, continuous lattices, Stone duality, split equalizer}


\pacs[MSC Classification]{06A11,06B35,06D22,18C15,18C20,18F60,18F70,54A99,54B30}

\maketitle

\section{Introduction}\label{sec1}

The importance and pervasiveness of the ideal lattice functor, or simply the ideal functor, in the study and history of Pointfree Topology are well-known, and its presence can be appreciated from the fact that it is to distributive lattices what the downset functor is to meet-semilattices (\citep[Remark 6.6]{BanPul2001}). In particular, it is a central tool in certain fundamental results: on one hand it provides a richer structure in the construction of compactifications of frames (\citep{Ban90,BezHar2014,FriScha2018}) and completions of uniform frames (\citep{Frith90,BanPul90}), and on the other it plays a pivotal role in the Stone representation of distributive lattices (\citep{Joh82}) as well as in the theory of dualities\footnote{``Just as in ring theory, the concepts of elements and ideals are of utmost importance in lattice theory.'' (\citep{HofMisStr74})} as exemplified by the dual adjunction that exists between distributive lattices and topological spaces via a schizophrenic object (\cite{Sim82,Joh82,HofMisStr74}). In general, functorial constructions from the ideal lattice produce either a monad or a coreflector, i.e. an idempotent comonad, in different parts of the existing literature. For instance, while P. Johnstone and H. Simmons consider the monad generated by the dual adjunction between spaces and lattices in \citep{Joh82} and \cite{Sim82} respectively, B. Banaschewski describes various reflective and coreflective subcategories with the help of the ideal lattice in \citep{Ban90}. These seemingly contradictory descriptions are not only a consequence of the change in the domain and codomain categories, but more precisely from a much deeper reality: they arise as successive iterations of the functor. This qualitative transformation between opposites can be understood as a dialectical process as explained by W. Lawvere in the examples provided in the paper \citep{Law96}. Here, the cylindrical model (\citep[Section I]{Law96}) of two sections united by a retraction (Cf. Lemma \ref{lemma: characterisation of algebras of lax idempotent monads}) qualitatively evolves into its opposite after another (quantitative) iteration of the functor.  We therefore concentrate our investigation on the monadic aspects of the ideal functor in relation to this iteration. Generally, a monad $\trip{T}=\mon{T}{m}{e}$ on a category \cat{C} induces an alternating chain of monads and comonads (\citep{Bar69}) on the various subcategories of algebras and coalgebras: 
$$
\dots \algn{T}{n}\to \coalgn{K}{n}\to \dots \to\algn{T}{1}\to\coalg{K}\to\alg{T}\to \cat{C}.
$$
The question naturally arises as to whether this sequence stops in a meaningful way or continues indefinitely. M. Barr has provided a complete answer in \citep{Bar69} for the category of sets, pointed sets and vector spaces over a field. In the category of distributive lattices endowed with the ideal functor, it actually follows from B. Jacobs' results (\cite[Theorem 4.5]{Bart2011}) on Kock-Z\"oberlein monads, also called lax idempotent monads, that \alg{T} and $\algn{T}{1}$ are equivalent.  

In this paper, we show that it is sufficient for \cat{C} to have either equalizers or coequalizers, for \alg{T} and $\algn{T}{1}$ to be equivalent in the sense that \alg{T} is monadic over \coalg{K}. Thus we reproduce B. Jacobs's result, but free of any condition on the monad \trip{T}. With an additional mild condition on \trip{T}, we can lift this equivalence to the level of \cat{C} and the category of free \trip{T}-algebras, hence generalising the Stone representation theorem for distributive lattices. 

The structure of the paper is very simple. We first introduce the necessary background in Section \ref{Preliminaries and background}, and then give a brief representation of algebras and coalgebras through split coequalizers and split equalizers in Section \ref{Presentation of algebras and coalgebras}. This representation is necessary to compensate for the absence of the chain of Galois connections used for lax idempotent monads (\citep{Hoff79,RoseWood2004,Bart94,Bart2011,Hofmann_Sousa2017}). The main results are then shown and discussed in Section \ref{Natural equivalence} together with a few illustrations. One key tool that we use is the Fakir construction of an idempotent approximation of a monad \citep{Fak70,LamRat73}. Its main idea stems from the observation that there is a natural and intimate relationship between \trip{T}-algebras and equalizers of pairs $(Te_X,e_{TX})$. The proofs then rely almost entirely on the manipulation of split equalizers and split coequalizers. In the case of the equivalence in Theorem \ref{theorem: equivalence between T-algebras and T1-algebras}, a split coequalizer is traced around the chain diagram described above. We reproduce here a step of B. Jacobs' proof of \cite[Theorem 4.5]{Bart2011}, and duly reference it, in order to properly separate it from the context of a lax idempotent monad and for the presentation to be self-contained. We end the paper with a note on projectives in the category \alg{T}; this is one categorical synthesis of the characterisations given in the papers \citep{Joh81,BanNie91,Dongsheng97,Mad91}.\\

\section{Preliminaries and background}
\label{Preliminaries and background}
\paragraph*{Distributive lattices and frames.}
A distributive lattice (\citep{Joh82,Sim82}) is a partially ordered set\footnote{In the context of the paper,``set'' means small set.} set admitting finite joins and finite meets\footnote{The lattice is then bounded.}, and for which the identity
$$
\meet{a}{(\join{b}{c})}=\join{(\meet{a}{b})}{(\meet{a}{c})}
$$
and its order-dual hold. A homomorphism between distributive lattices is understood as a function that preserves the operations $\wedge$ and $\vee$, including the top element $\wedge\emptyset=1$ and the bottom element $\vee\emptyset=0$. These form the category \cat{DLat}.

A frame (\citep{PicPul2012,Joh82}) $L$ is a complete lattice with the equational identity
$$
\meet{a}{\left(\ds{\bigvee} S\right)}=\ds{\bigvee\set{\meet{a}{s}\ |\ s\in S}}
$$
for all $a\in L$ and any $S\subseteq L$. Frame homomorphisms or frame maps are functions \map{f}{L}{M} that preserve finite meets  and arbitrary joins. Frames and their homomorphisms form a category denoted by \cat{Frm}. There is forgetful functor $\cat{Frm}\to\cat{DLat}$.\\
\begin{definition}
A subset $J\subseteq L$ is called an ideal if $J$ is a downset, i.e. $b\in J$ provided there is $a\in J$ with $b\leq a$, and if it is closed under the formation of finite joins. 
\end{definition}
\noindent
The set of all ideals on a frame $L$ is denoted by $\id L$. $\id L$ is a frame with the following operations: meets are given by set-intersections and joins given by
$$
\bigvee\mathscr{J}=\bigcup\set{I_1\vee I_2\vee\dots\vee I_n\ |\ I_1, I_2, \dots, I_n\in\mathscr{J}\text{ and }n\in \mathbb{N}},
$$
where $I_1\vee I_2\vee\dots\vee I_n=\set{i_1\vee i_2\vee\dots\vee i_n\ |\ i_k\in I_k\text{ for }1\leq k\leq n}$. Thus, if $\mathscr{J}$ is a directed set, then $\ds{\bigvee}\mathscr{J}=\bigcup\mathscr{J}$. If \map{f}{L}{M} is a  homomorphism between distributive lattices, then the function \map{\id f}{\id L}{\id M} given by
$$
\id f(I)=\set{b\ |\ b\leq f(a)\text{ for some }a\in I},
$$
is a frame homomorphism. In particular \id\ is an endofunctor on \cat{DLat}. The frames of the form $\id D$, where $D$ is a distributive lattice, are called {\em coherent frames} (\citep{Ban81,Joh82}). 
\noindent
We are interested in the way-below relation $\ll$ which is defined by 
$$
a\ll b\text{ if for any }S\subseteq L,\text{ if }b\leq\bigvee S\text{ then } a\leq\bigvee F\text{ for some finite }F\subseteq S.
$$ 
A frame $L$ is said to be {\em stably compact} (\citep{Sim82,BanBru,Jun2004}) if $\ll$ is a sublattice of $L\times L$ and if for all $a\in L$, $a=\bigvee\set{b\ |\ b\ll a}$. The category of stably compact frames and frame homomorphisms that preserve\footnote{These are called proper maps or perfect maps. The condition is equivalent to the right adjoint of the frame map preserving directed joins.} $\ll$ is denoted by \cat{StKFrm}. The category of coherent frames and frame homomorphisms preserving {\em compact elements}, i.e. elements satisfying $a\ll a$, is denoted by \cat{CohFrm}. We have inclusion functors $\cat{CohFrm}\to\cat{StKFrm}\to\cat{Frm}$.

\paragraph*{Notations for natural transformations.} Given two pairs of functors \map{F,G}{\cat{C}}{\cat{D}} and \map{H,K}{\cat{D}}{\cat{E}} and two natural transformations $\alpha:F\to G$ and $\beta:H\to K$, the horizontal composition $\beta\circ\alpha:HF\to KG$ is defined by 
$$
(\beta\circ\alpha)_X=\beta_{GX}H(\alpha_X)=K(\alpha_X)\beta_{FX}
$$
for each $X$ in \cat{C}, for which we simply write $\beta\circ\alpha=\beta G\cdot H\alpha=K\alpha\cdot\beta F$.
If \map{L}{\cat{C}}{\cat{D}} is another functor and $\delta:G\to L$ another natural transformation, then the vertical composition $\comp{\delta}{\alpha}:F\to L$ is simply defined as $\comp{\delta_X}{\alpha}_X$ for each $X\in\cat{C}$.\\

\begin{lemma}
\label{lem: middle-interchange law}
(Middle-Interchange Law) Given natural transformations
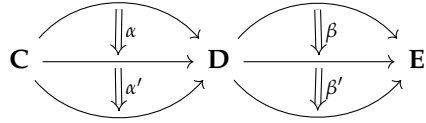
\[
\begin{tikzcd}[column sep=2cm]
 \cat{C}   \arrowtcmid{r}{} \arrowtcupp{r}{\alpha}\arrowtclow{r}{\alpha\smash'} & \cat{D} \arrowtcmid{r}{} \arrowtcupp{r}{\beta}\arrowtclow{r}{\beta\smash'} & \cat{E}
\end{tikzcd}\]
we have $(\beta\smash'\circ\alpha\smash')\cdot(\beta\circ\alpha)=(\beta\smash'\cdot\beta)\circ(\alpha\smash'\cdot\alpha)$. (\cite{Low,Mac}.)
\end{lemma}

\paragraph*{Monads.} Recall that a {\em monad} (\citep{,EilMoo65,Low,Mac,Rie}) $\mathbb{T}$ on a category \cat{C} is a triple \mon{T}{m}{e}, where \map{m}{TT}{T} and \map{e}{1}{T} are natural transformations satisfying the identities 
$$
\comp{m}{Tm}=\comp{m}{m T}\text{ and }\comp{m}{e T}=\comp{m}{Te}=1_T.
$$
\begin{example}
\label{example: ideal functor as a monad}
The ideal functor \id, together with the natural transformations provided by the set-union\footnote{An ideal on \id L is by definition directed.} \map{\bigcup}{\id\id}{\id} and the principal ideal map \map{\downarrow}{1}{\id}, is a monad (Cf. \citep[Section 2]{Sim82} and \citep[Exercise 4.6]{Joh82}).\\
\end{example}
A \trip{T}-algebra (or an {\em Eilenberg-Moore algebra}) is a pair $(X,a)$, where $X\in\cat{C}$ and \map{a}{TX}{X} a morphism such that 
$$
\comp{a}{Ta}=\comp{a}{m_X}\text{ and }\comp{a}{e_X}=1_X.
$$
If $(X,a)$ and $(Y,b)$ are \trip{T}-algebras, then a \trip{T}-algebra homomorphism \map{f}{(X,a)}{(Y,b)} is a morphism \map{f}{X}{Y} in \cat{C} such that $\comp{f}{a}=b\cdot Tf$. The category of \trip{T}-algebras and \trip{T}-algebra homomorphisms is denoted by \powcat{\cat{C}}{T} or by \alg{T}. The forgetful functor $\map{\powcat{G}{T}}{\powcat{\cat{C}}{T}}{\cat{C}}:(X,a)\mapsto X$ admits a left adjoint $\map{\powcat{F}{T}}{\cat{C}}{\powcat{\cat{C}}{T}}$ defined by $\powcat{F}{T}(\map{f}{X}{Y})=\map{Tf}{(TX,m_X)}{(TY,m_Y)}$. The unit of this adjunction is given by $e_X:X\to \powcat{G}{T}\powcat{F}{T}X$ and the co-unit is provided by the algebra morphisms $\varepsilon_{(X,a)}=a:\powcat{F}{T}\powcat{G}{T}(X,a)\to(X,a)$. It is important to note in addition that \powcat{G}{T} reflects isomorphisms.\\
\begin{definition}
\label{definition: split coequalizer}
(\citep{Low,Rie}) The diagram
$$\xymatrix{
A \ar@<0.5ex>[r]^{f} \ar@<-0.5ex>[r]_{g} & B\ar@/^1.4pc/[l]^{t} \ar[r]^q & C\ar@/^1pc/[l]^{s}.
}$$
is said to be a split coequalizer if $q\cdot f=q\cdot g$, $q\cdot s=1_C$, $f\cdot t=1_B$ and $s\cdot q=g\cdot t$.\\
\end{definition}
\noindent
Split coequalizers are absolute, i.e. they are preserved by any functor. 

Any adjunction $(F\dashv G,\eta,\epsilon)$ where \map{G}{\cat{A}}{\cat{C}} and \map{F}{\cat{C}}{\cat{A}} induces a monad $(GF,G\epsilon F,\eta)$ on \cat{C} (\citep{EilMoo65}). If the adjunction $F\dashv G$ induces the same monad \trip{T} on \cat{C}, that is $GF=T, G\epsilon F=m$ and $\eta=e$, then there is a unique functor \map{C}{\cat{A}}{\powcat{\cat{C}}{T}} - called {\em comparison functor}, such that $CF=\powcat{F}{T}$ and $\powcat{G}{T}C=G$. We say that \cat{A} is monadic (resp. strictly monadic) over \cat{C} if $C$ is an equivalence (resp. isomorphism). Beck's monadicity criterion (\citep{Beck2003,Low,Rie}) implies in this case that $G$ creates coequalizers of $G$-split pairs: i.e. if $f=Gu$ and $g=Gv$ in Definition \ref{definition: split coequalizer}, then $u$ and $v$ admit a coequalizer $r$ in \alg{T} such that the $G$-image of the diagram $\comp{r}{u}=\comp{r}{v}$ is isomorphic to the split coequalizer. This relation of isomorphism is strict for \powcat{G}{T}.

The adjunction $\powcat{F}{T}\dashv\powcat{G}{T}$ induces a comonad $\trip{K}=(T,\powcat{F}{T}e\powcat{G}{T},\varepsilon)$ on the category of algebras \alg{T}, where $\powcat{F}{T}e\powcat{G}{T}=\map{Te}{(T,m)}{(TT,mT)}$. This generates an alternating sequence of monads and comonads as follows:

$$\xymatrix{
\cat{C}\ar@(ul,ur)^{T} \ar@<1.0ex>[dd]^{\powcat{F}{T}}_*-<0.55ex>{_\vdash} & & \cat{Alg}(\mathbb{T}_1)\ar@(ul,ur)^{K_1} \ar@<1.0ex>[dd]^{\ds{G}^{\mathbb{T}_1}}_*-<0.55ex>{_\dashv} \ar@<-1.0ex>[rr]_{\ds{F}^{\mathbb{K}_1}}^*-<0.55ex>{_\perp} & & \cat{Coalg}(\mathbb{K}_1)\ar@(ul,ur)^{T_2}  \ar@<-1.0ex>[ll]_{\ds{G}^{\mathbb{K}_1}} \ar@{.>}@<1.0ex>[dd]\\ \\
\alg{T}\ar@(dl,dr)_{K} \ar@<1.0ex>[uu]^{\powcat{G}{T}} \ar@<-1.0ex>[rr]_{\powcat{F}{K}}^*-<0.55ex>{_\perp} & & \cat{Coalg}(\trip{K})\ar@(dl,dr)_{T_1}  \ar@<1.0ex>[uu]^{\ds{F}^{\mathbb{T}_1}} \ar@<-1.0ex>[ll]_{\powcat{G}{K}} & & \dots \ar@{.>}@<1.0ex>[uu]
}$$
We have sequences of comparison functors \map{\powcat{F_{n}}{T}}{\algn{T}{n-1}}{\algn{T}{n}} and \map{\powcat{C_{n}}{K}}{\coalgn{K}{n-1}}{\coalgn{K}{n}} for $n\geq 1$, where $\trip{T}_0=\trip{T}$ and $\trip{K}_0=\trip{K}$.\\

\begin{definition}
A category \cat{C} is said to be order-enriched if each hom-sets $\text{Hom}(X;Y)$ is a partially ordered set for any pair $X,Y\in\cat{C}$ and the composition of morphisms preserves the order. The monad \trip{T} is said to be order-enriched if the underlying functor $T$ preserves the order.\\
\end{definition}
\begin{definition}
(\citep{Low,Koc}) An order-enriched monad \trip{T} is said to be lax idempotent or a Kock-Z\"oberlein monad if $Te_X\leq e_{TX}$ for any $X\in\cat{C}$.\\
\end{definition}
\noindent
We will omit mentioning that \trip{T} is order-enriched each time we consider a lax idempotent monad, as this is clear from the context\footnote{Unlike in \citep{Low}, we consider the order to be separated, i.e. $x\leq y$ and $y\leq x$ imply $x=y$.}.\\
\begin{lemma}
\label{lemma: characterisation of algebras of lax idempotent monads}
(\citep[Section II.4.9]{Low}) For a monad \trip{T}, we have the following equivalences for any $X\in\cat{C}$:
$$
Te_X\leq e_{TX}\Longleftrightarrow Te_X \dashv m_X\Longleftrightarrow m_X\dashv e_{TX}.
$$
In such a case, any morphism \map{a}{TX}{X} with $\comp{a}{e_X}=1_X$ necessarily satisfies $a\dashv e_X$ and this makes $(X,a)$ a \trip{T}-algebra.\\
\end{lemma}
\noindent
We mention the following contribution by B. Jacobs.\\
\begin{theorem}
\label{theorem: Bart Jacobs' main result}
(\citep{Bart2011}) For a lax idempotent monad \trip{T} on a category \cat{C} with equalisers, one has $\alg{T}\simeq\algn{T}{1}$.
\end{theorem}

\section{Presentation of algebras and coalgebras}
\label{Presentation of algebras and coalgebras}
A \trip{T}-algebra $(X,a)$ can be presented as in the following split coequalizer (\citep{Low,Rie}):
$$\xymatrix{
TTX \ar@<0.5ex>[r]^{Ta} \ar@<-0.5ex>[r]_{m_X} & TX\ar@/^1.4pc/[l]^{eT} \ar[r]^a & X\ar@/^1pc/[l]^{e}.
}$$
A \trip{K}-coalgebra structure on $(X,a)$ is then a \trip{T}-algebra morphism \map{c}{(X,a)}{(TX,m_X)} that is part of the following split equalizer:
$$\xymatrix{
(X,a)\ar[r]^{c} & (TX,m_X)\ar@<0.5ex>[r]^{Tc} \ar@<-0.5ex>[r]_{Te}\ar@/^1pc/[l]^{a}  & (TTX,m_{TX})\ar@/^1.4pc/[l]^{m_X}.
}$$ 
We shall denote this \trip{K}-coalgebra as an ordered triple $(X,a,c)$. Thus a $\trip{T}_1$-algebra structure on $(X,a,c)$ is a \trip{K}-coalgebra morphism \map{b}{(TX,m_X,Te_X)}{(X,a,c)} that is a split coequalizer:
$$\xymatrix{
(TTX,m_{TX},Te_{TX}) \ar@<0.5ex>[r]^{Tb} \ar@<-0.5ex>[r]_{Ta} & (TX,m_X,Te_X)\ar@/^1.4pc/[l]^{Te} \ar[r]^b & (X,a,c)\ar@/^1pc/[l]^{c}. 
}$$
A $\trip{T}_1$-algebra as described above is denoted by $(X,a,c,b)$. Thus free \trip{K}-algebras and free $\trip{T}_1$-algebras are given by $(TX,m_X,Te_X)$ and $(TX,m_X,Te_X,Ta)$ respectively so that the functor \powcat{F}{T} restricts itself to the comparison functor \map{\powcat{F_{1}}{T}}{\alg{T}}{\algn{T}{1}} between the adjunctions that induce the monad $\trip{T}_1$. It is a restriction in a sense that 
$$
(\powcat{G}{K}G^{\trip{T}_1})\cdot\powcat{F_{1}}{T} = K = \powcat{F}{T}\cdot\powcat{G}{T}.
$$
For a lax idempotent monad, the description of the algebras and coalgebras can be greatly simplified. In the table below, we list the identities characterising the algebras and coalgebras for a lax idempotent monad.\\
\begin{center}
\begin{tabular}{ |p{5cm}|p{5cm}| }
 \hline
 \textbf{(co)Monads } & \textbf{(co)Algebras}\\
 \hline
  $\trip{T}=(T,m,e)$ &  $a\dashv e_X$ and $a\cdot e_X=1$\\ 
 \hline
 $\trip{K}=(T,Te,a)$ &   $c\dashv a$ and $a\cdot c=1$ \\ 
 \hline
 $\trip{T}_1=(T,Ta,c)$ &   $b\dashv c$ and $b\cdot c=1$ \\ 
 \hline
 $\cdots$ & $\cdots$ \\
\end{tabular}
\end{center} 
The simplification in the table is due to various work that can be traced back to the work of R.-E. Hoffman on continuous posets \citep{Hoff79}, the papers by A. Kock (\citep{Koc}) and V. Z\"oberlein (\cite{Zob76}) on lax idempotent monads. These characterisations are also described by B. Jacobs in \citep{Bart2011}.\\

\begin{remark}
\label{remark: on split algebras}
For a lax idempotent monad \trip{T}, the fact that \map{c}{X}{TX} is as well a \trip{T}-algebra morphism is due to its being a left adjoint of $a$. This is given by \citep[Proposition 2.5]{Koc} and \citep[Proposition 3.8]{Hofmann_Sousa2017}. Split algebras (Cf. \citep[Definition 3.6]{Hofmann_Sousa2017} and \citep{RoseWood2004}) are then characterised as coalgebras of the comonad \trip{K}.\\
\end{remark}
\noindent
We mention the following:\\
\begin{theorem}
\label{theorem: characterisation of algebras and coalgbras for the ideal monad}
\begin{enumerate}
\item A distributive lattice $D$ is a frame if and only if \map{\downarrow}{D}{\id D} admits a left adjoint $a$ such that $\comp{a}{\downarrow}=1_D$. In this case, the left adjoint is given by the join $\ds{\bigvee}$.
\item A frame $L$ is stably compact if and only if the join \map{\ds{\bigvee}}{\id L}{L} admits a left adjoint $c$ such that $\comp{\ds{\bigvee}}{c}=1_D$. Here $c(x)=\set{y\in L\ |\ y\ll x}$\footnote{This expression is forced by the identities $\comp{\ds{\bigvee}}{c}=1_D$ and $\comp{c}{\ds{\bigvee}}\leq 1_{\mathfrak{I}D}$.}.\\
\end{enumerate}
\end{theorem}
\noindent
With the help of Lemma \ref{lemma: characterisation of algebras of lax idempotent monads}, the above characterisation can be deduced from the facts that \cat{Frm} is the category of algebras for the ideal monad (\citep[Exercise 4.6]{Joh82}), and that \cat{StKFrm} is the category of coalgebras of the comonad \mon{\id}{\id(\downarrow)}{\ds{\bigvee}} as shown in \citep[Section 3]{BanBru}. It is known since \citep[Theorem 4.2]{Koc} and \citep{Hoff79} that the way below relation $\ll$ on a poset is directly linked to the condition of a certain join map admitting a left adjoint. We refer the reader to the compendium \citep{Gie}, as well as the seminal paper \citep{Scott72} by D. Scott and the paper  \citep{Day75} by Alan Day for a brief backgroound on continuous lattices. We note that D. Scott was the first to introduce and systematically study continuous lattices. The category of algebras of the induced monad \mon{\id}{\id(\bigvee)}{c} is also familiar to us:\\

\begin{theorem}
\label{theorem: characterisation of coherent frames as algebras}
The algebras of the monad \mon{\id}{\id(\bigvee)}{c} are characterised by those stably compact $L$ for which \map{c_L}{L}{\id L} admits a left adjoint \map{b}{\id L}{L} such that $\comp{b}{c_L}=1_L$.
\end{theorem}
\begin{proof}
This can be obtained by following an argument similar to \citep[Theorem 3.4, Lemma 3.6]{Bart94} and \citep[Theorem 2.6]{Hoff79}.
\end{proof}
\noindent
We observe that $b$, being a left adjoint to $c$, is given by:
$$ b(J)=\ds{\bigvee}\set{y\in J\ |\ y\ll y} .$$
$L$ then becomes a coherent frame, that is $L=\id D$ where $D=\set{y\in L\ |\ y\ll y}$. However, not all coherent frames satisfy Theorem \ref{theorem: characterisation of coherent frames as algebras}. The following result is at the heart of the equivalence in Theorem \ref{theorem: Bart Jacobs' main result}. \\
\begin{proposition}
\label{proposition: characterisation of T1-algebras for a lax idempotent monad}
(\citep{Bart2011}) Let \trip{T} be a lax idempotent monad and let us assume that $e_X$ is the equalizer of the pair $(e_{TX},Te_X)$ for a fixed $X$. It follows that $X$ admits a \trip{T}-algebra structure if and only if $TX$ admits a $\trip{T}_1$-algebra structure.
\end{proposition}
\begin{proof}
The necessary condition is obvious. Conversely, suppose that $b\cdot Te_X=1_{TX}$ exists. As shown in \citep[Theorem 4.5]{Bart2011}, the sequence $b\dashv Te_X\dashv m_X\dashv e_{TX}$ makes $e_X$ a split monomorphism: we have
\begin{align*}
\comp{\comp{Te_X}{b}}{e_{TX}} &= \comp{\comp{Tb}{Te_{TX}}}{e_{TX}}\text{ (}b\text{ is a coalgebra map)}\\
&= \comp{\comp{Tb}{e_{TTX}}}{e_{TX}}\\
&= \comp{\comp{e_{TX}}{b}}{e_{TX}}.
\end{align*}
There is a unique morphism \map{r}{TX}{X} such that $\comp{e_X}{r}=\comp{b}{e_{TX}}$. Composing with $e_X$ from the  right gives $e_X\cdot\comp{r}{e_X}=\comp{\comp{b}{Te_X}}{e_X}=e_X$. One has $\comp{r}{e_X}=1$ (\citep{Bart2011}).
\end{proof}
\begin{remark}
\label{remark: downset map as an equalizer}
For a distributive lattice $D$, since \cat{DLat} is monadic over \cat{Set} (\citep[Proposition 3.7]{Joh82}), \map{\downarrow}{D}{\id D} is an equalizer of $(\downarrow_{\mathfrak{I}D},\id(\downarrow_D))$. Consequently, a coherent frame $\id D$ is an algebra of the monad \mon{\id}{\id(\bigvee)}{c} if and only if $D$ is a frame. A similar situation is described in \citep[Corollary 2.7]{Hoff79} for algebraic lattices, and subsequently discussed in \citep[Section 8]{RoseWood2004}.
\end{remark}

\section{Natural equivalences}
\label{Natural equivalence}
In order to show that under some mild conditions \powcat{F_1}{T} is an equivalence, we shall introduce a functor \map{\Phi}{\powcat{F}{T}(\cat{C})}{\cat{C}} that is provided by the Fakir construction of idempotent approximation of a monad (\cite{Fak70}). For a given monad $\trip{T}=\mon{T}{m}{e}$, consider the equalizer 
$$
\eq{\fak{T}}{\varphi}{T}{Te}{eT}{TT}.
$$
\begin{theorem}
\label{theorem: Fakir construction}
(Fakir \citep{Fak70}) The functor \fak{T} underlies a monad $\trip{I}_T=\mon{\fak{T}}{\fak{m}}{\fak{e}}$ where \fak{m} is unique such that $\comp{\varphi}{\fak{m}}=\comp{m}{(\ocomp{\varphi}{\varphi})}$ and \fak{e} is unique such that $\comp{\varphi}{\fak{e}}=e$.\\
\end{theorem}
\noindent
We mention the following results which are of interest to us.\\
\begin{proposition}
\label{proposition: results from Fakir construction}
(\citep[Proposition 3]{Fak70}) \map{\fak{e}T}{T}{\fak{T}T} is an isomorphism, and the following are equivalent:
\begin{enumerate}
\item[•] $\trip{I}_T$ is idempotent;
\item[•] \map{T\fak{e}}{T}{T\fak{T}} is an isomorphism;
\end{enumerate}
\end{proposition}
Now consider the category $\powcat{F}{T}(\cat{C})$ of free \trip{T}-algebras and morphisms $\powcat{F}{T}f$ where $f$ is in \cat{C}. We consider the functor \map{\powcat{\theta}{T}}{\cat{C}}{\powcat{F}{T}(\cat{C})} so that if $\powcat{E}{T}$ is the inclusion functor  $\powcat{F}{T}(\cat{C})\to\alg{T}$, then $\comp{\powcat{E}{T}}{\powcat{\theta}{T}}=\powcat{F}{T}$. 
\noindent
In a similar manner, we consider \map{\powcat{\theta_1}{T}}{\alg{T}}{\powcat{F_1}{T}(\alg{T})}, where \powcat{F_1}{T}(\alg{T}) is the category of free $\trip{T}_1$-algebras with morphisms of the form $\powcat{F_1}{T}f$, as well as the inclusion $\powcat{E_1}{T}:\powcat{F}{T}(\cat{C})\to\alg{T}$. We have 
$$\xymatrix{
\cat{C}\ar[rr]^{\theta^{\mathbb{T}}} & & F^{\mathbb{T}}(\cat{C}) \ar[rr]^{E^\mathbb{T}} & & \alg{T} \\
\alg{T}\ar[rr]^{\theta_1^{\mathbb{T}}} \ar[u] & & F_1^{\mathbb{T}}(\alg{T}) \ar[rr]^{E_1^\mathbb{T}} \ar[u] & & \algn{T}{1} \ar[u]
}$$
where vertical functors are forgetful. \\
\begin{definition}
\label{definition: Fakir basis}
We define \map{\Phi}{\powcat{F}{T}(\cat{C})}{\cat{C}} by $\Phi(\powcat{F}{T}f)=\fak{T}f$ for all $f$ in \cat{C}. In particular $\Phi\cdot\powcat{\theta}{T}=\fak{T}$. This gives two natural transformations \map{\fak{e}}{1}{\Phi\cdot\powcat{\theta}{T}} and \map{\powcat{\theta}{T}\fak{e}}{1}{\comp{\powcat{\theta}{T}}{\Phi}}.\\
\end{definition}

\begin{theorem}
\label{theorem: monadic form of stone-type dualites}
If $\trip{I}_T$ is an idempotent monad, then $\powcat{\theta}{T}(\fak{e})$ is a natural isomorphism. In particular if $\trip{I}_T$ is the identity monad, then we have an equivalence between the category of free \trip{T}-algebras and the ambient category \cat{C}.
\end{theorem}
\begin{proof}
We note that $\powcat{\theta}{T}\fak{e}=T\fak{e}:(T,m)\to(T\fak{T},m\fak{T})$. Since \powcat{G}{T} reflects isomorphisms, $T\fak{e}$ is an isomorphism if and only if \powcat{\theta}{T}\fak{e} is. The result easily follows form Proposition \ref{proposition: results from Fakir construction}. 
\end{proof}
\noindent
We shall restrict $\Phi$ to $F_1^{\mathbb{T}}(\alg{T})\to\alg{T}$. We first give a strengthening of the fact that $\fak{e}T$ is an isomorphism.\\
\begin{lemma}
\label{lemma: Fakir construction fixes T-algebras}
For a \trip{T}-algebra $(X,a)$, $e^{\varphi}_X$ is an isomorphism. 
\end{lemma}
\begin{proof}
If $(X,a)$ is a \trip{T}-algebra, then we have a split equalizer
$$\xymatrix{
X\ar[r]^{e} & TX\ar@<0.5ex>[r]^{Te} \ar@<-0.5ex>[r]_{eT}\ar@/^1pc/[l]^{a}  & TTX\ar@/^1.4pc/[l]^{Ta},
}$$ 
which makes $e^{\varphi}_X$ an isomorphism. 
\end{proof}
\noindent
This allows us to define the restriction \map{\Phi_1}{F_1^{\mathbb{T}}(\alg{T})}{\alg{T}} of $\Phi$ by $\Phi_1(TX,m_X,Te,Ta)=(X,a)$ and $\Phi_1(\powcat{F_1}{T}f)=f$.\\
\begin{lemma}
\label{lemma: equivalence between split equalizers and split equalizers}
Suppose that \cat{C} has either equalizers or coequalizers. If $(X,a,c,b)$ is a $\trip{T}_1$-algebra, then we have the following split coequalizer and split equalizer
 $$\xymatrix{
TX \ar@<0.5ex>[r]^{b} \ar@<-0.5ex>[r]_{a} & X\ar@/^1.4pc/[l]^{e} \ar[r]^{r} & X_c\ar@/^1pc/[l]^{\kappa} \ar[r]^{\kappa} & X\ar@<0.5ex>[r]^{c} \ar@<-0.5ex>[r]_{e}\ar@/^1pc/[l]^{r}  & TX\ar@/^1.4pc/[l]^{b}
}$$
\end{lemma}
\begin{proof}
The identities $\comp{a}{e}=1$ and $\comp{b}{c}=1$ hold as $(X,a,b,c)$ is a $\trip{T}_1$-algebra. Suppose that \cat{C} has equalizers\footnote{The case where it has coequalizers is treated analogously.} and denote by $\kappa$ the equalizer of $(c,e)$. As shown in the proof of \citep[Theorem 4.5]{Bart2011}: $c\cdot b\cdot e=\comp{\comp{Tb}{Te}}{e}=\comp{\comp{Tb}{eT}}{e}=\comp{\comp{e}{b}}{e}$; there is a unique $r$ such that $\kappa\cdot r=b\cdot e$. One shows that $\comp{\comp{k}{r}}{k}=\comp{\comp{b}{c}}{k}=k$ and so $r\cdot \kappa=1$. We have a split equalizer on the right. It remains to show that $r\cdot a=r\cdot b$ to have a split coequalizer on the left. We have
\begin{align*}
\comp{\comp{\kappa}{r}}{b} &= \comp{\comp{b}{e}}{b}\\
&= \comp{\comp{b}{Tb}}{eT}\\
&= \comp{\comp{b}{Ta}}{eT}\text{ (}b\text{ coequalizes }Ta\text{ and }Tb\text{)}\\
&= \comp{\comp{b}{e}}{a}\\
&= \comp{\comp{\kappa}{r}}{a}.
\end{align*}
Since $\kappa$ is a monomorphism, the identity $r\cdot a=r\cdot b$ follows.
\end{proof}

\begin{theorem}
\label{theorem: equivalence between T-algebras and T1-algebras}
If \cat{C} admits either coequalizers for \trip{T}-algebra morphisms or equalizers for any parellel pair of morphisms, then \powcat{F_1}{T} is an equivalence.
\end{theorem}
\begin{proof}
First, since the restriction of $\trip{I}_T$ on \alg{T} is essentially the identity monad, we have $\comp{\Phi_1}{\powcat{\theta_1}{T}}\cong 1$ and $\comp{\powcat{\theta_1}{T}}{\Phi_1}\cong 1$. Now, let $(X,a,c,b)$ be a $\trip{T}_1$-algebra and consider the split coequalizer from Lemma \ref{lemma: equivalence between split equalizers and split equalizers}:
 $$\xymatrix{
TX \ar@<0.5ex>[r]^{b} \ar@<-0.5ex>[r]_{a} & X\ar@/^1.4pc/[l]^{e} \ar[r]^r & X_c\ar@/^1pc/[l]^{\kappa}
}$$
As an absolute coequalizer, this is preserved by the comparison functor from \cat{C} to \coalg{K} which arises from the adjunctions inducing the comonad \trip{K} on \alg{T}. Thus we have a split coequalizer
$$\xymatrix{
(TTX,m_{TX},Te_{TX}) \ar@<0.5ex>[r]^{Tb} \ar@<-0.5ex>[r]_{Ta} & (TX,m_X,Te_X)\ar@/^1.4pc/[l]^{Te} \ar[r]^{Tr} & (TX_c,m_{X_c},Te_c)\ar@/^1pc/[l]^{T\kappa}. 
}$$
in \coalg{K}. Now, both $Ta$ and $Tb$ are morphisms of (free) $\trip{T}_1$-algebras. By Beck's monadicity theorem, since \algn{T}{1} is monadic over \coalg{K} in the strict sense, there is a $\trip{T}_1$-algebra structure $\gamma$ on $(TX_c,m_{X_c},Te_c)$ and $Tr$ becomes a $\trip{T}_1$-algebra morphism. The universality of the coequalizer $b$ in \algn{T}{1} shows that
$$
(TX_c,m_{X_c},Te_c,\gamma)\cong (X,a,c,b)\text{ and } Tr\cong b\cong\gamma.
$$
Finally since \alg{T} is monadic over \cat{C}, $X_c$ admits a \trip{T}-algebra structure $a_c$, as both $a$ and $b$ are \trip{T}-algebra morphisms. The isomorphism $TX_c\cong X$ implies that $r\cong a_c$ as morphisms of \trip{T}-algebras. If \map{f}{(X,a,c,b)}{(X',a',c',b')} is a $\trip{T}_1$-algebra morphism, then we form the (split) equalizers $\kappa,\kappa'$ and (split) coequalizers $r,r'$. There is a unique morphism $f_c$ such that $\comp{f_c}{r}=\comp{r'}{f}$ and $\comp{f}{\kappa}=\comp{\kappa'}{f_c}$. Thus $f_c$ is a \trip{T}-algebra morphism and $\theta_1^{\mathbb{T}}f_c\cong f$. It follows that the inclusion \powcat{E_1}{T} is an equivalence.
\end{proof}

\begin{corollary}
\label{corollary: equivalence between T-algebras and T1-algebras}
With the assumption of Theorem \ref{theorem: equivalence between T-algebras and T1-algebras}, \powcat{F}{K} is monadic, $F^{\mathbb{T}_1}$ is comonadic and the action of $\trip{K}_1$ on \algn{T}{1} is equivalent to that of \trip{K} on \alg{T}. \\
\end{corollary}

We note that B. Jacobs' result in \citep{Bart2011} requires that certain equalizers exist in \cat{C} in order to construct the object $X_c$. The identification of the algebra $(X,a,c,b)$ with $(TX_c,m_{X_c},Te_c,Ta_c)$ proceeds from the fact that for a lax idempotent monad, the adjoints $a\dashv c\dashv b$ coincide with $m_{X_c}\dashv Te_c\dashv Ta$ due to the property of universality of Galois connections. The verification of the equivalence in Theorem \ref{theorem: Bart Jacobs' main result} then requires one to check that the relevant diagrams commute. The absence of the Galois connections has been remedied here by the use of the characterisation of algebras through split coequalizers, reaffirming their importance in the theory of monads.\\

\begin{example}
\label{example: ideal functor and downset functor}
We are mostly interested in the ideal functor and the down set functor. Similar examples that are treated by Bart Jacobs in \citep{Bart94} and \citep{Bart2011} are assumed to be known. These functors have the common property that they are lax idempotent, and so the results in  \citep{Bart94} apply. Here, we augment these results with the  equivalence from Theorem \ref{theorem: monadic form of stone-type dualites}.
\begin{enumerate}
\item[(i)] In the category \cat{DLat}, the first inductive step \mon{\fak{\id}}{\fak{\bigvee}}{\fak{\downarrow}} of Fakir construction  is essentially the identity monad (Cf. Remark \ref{remark: downset map as an equalizer}). The equivalence between \cat{DLat} and the category \cat{CohFrm} then follows from Theorem \ref{theorem: monadic form of stone-type dualites}. This is part of the Stone representation for distributive lattices (\citep[Section 3.3]{Joh82}, \citep{Ban81}). We emphasize the fact that \cat{DLat} is complete and so Theorem \ref{theorem: equivalence between T-algebras and T1-algebras} applies as well. A case that is analogous to this is that of the category $\sigma\cat{Frm}$ (\citep{BanMat2003}) of $\sigma$-frames with the monad $\mathscr{H}$ taking a $\sigma$-frame $L$ to the frame $\mathscr{H}L$ of its $\sigma$-ideals . This establishes an equivalence between $\sigma\cat{Frm}$ and the category $\sigma\cat{CohFrm}$ of $\sigma$-coherent frames via Theorem \ref{theorem: monadic form of stone-type dualites}. The functor $\Phi$ (Definition \ref{definition: Fakir basis}) is given by the functor $\mathrm{Lind}$ that extracts the Lindel\"{o}f elements from the free algebras $\mathscr{H}L$.\\

\item[(ii)] Consider the category \cat{MLat} of meet-semilattices and meet-semilattice homomorphisms with the downset functor \map{\mathfrak{D}}{\cat{MLat}}{\cat{MLat}}. This form a monad \trip{T} with the natural transformations given by set-theoretic unions \map{\bigcup}{\mathcal{D}\mathcal{D}}{\mathcal{D}} and downset maps \map{\downarrow}{1}{\mathcal{D}}. The category of algebras \alg{T} is \cat{Frm} and the category \coalg{K} is the category of supercontinuous frames with their corresponding homomorphisms (\citep{BanNie91,BanPul2001}). Here \mon{\fak{\mathfrak{D}}}{\fak{\bigvee}}{\fak{\downarrow}} is essentially the identity. Theorem \ref{theorem: monadic form of stone-type dualites} shows that \cat{MLat} is equivalent to the category \cat{SCohFrm} of supercontinuous frames (Cf. \citep[Remark 3]{BanNie91} and \citep[Proposition 6.2]{BanPul2001}).\\

\item[(iii)] In \citep{Dongsheng97}, Z. Dhongsheng considers the functor \map{Z}{\cat{MLat}}{\cat{MLat}} which is a submonad \trip{T} of $\mathfrak{D}$ from the previous example\footnote{That this is a submonad is straightforward, given the requirements of Definition 1.1 in \citep{Dongsheng97}}. From this follow three subcategories: the category of $Z$-frames, the category of stably $Z$-continuous semilattices, and the category of coherent $Z$-frames. These are respectively the category \alg{T}, the category \coalg{K} and the category of free \trip{T}- algebras. Theorem \ref{theorem: monadic form of stone-type dualites} applies here as well.\\
\end{enumerate}
\end{example}

\begin{example}
\label{example: ultarfilter and filter monad functors}
The following classical examples trivially satisfy the condition of Theorem \ref{theorem: equivalence between T-algebras and T1-algebras}: the category of topological spaces and continuous maps with the ultrafilter monad or with the prime open filter monad ((\citep{Raz2022,Sim82})), the category of $T_0$-spaces and continuous maps with the open filter monad (\citep{Hofmann_Sousa2017,Esc98,Esc2002}). For those which are lax idempotent, B. Jacobs' result applies to these instances. Theorem \ref{theorem: equivalence between T-algebras and T1-algebras} allows for the extension of the results to arbitrary monads. However, for those classical instances we do not know if the first inductive step of the Fakir idempotent approximation reduces to the identity monad.\\
\end{example}

\begin{example}
\label{example: Michale Barr}
M. Barr considers the following examples in \citep{Bar69}: \cat{Set}, the category $(1,\cat{Set})$ of pointed sets, and the category $\cat{V}_K$ of vector spaces over a field $K$ and looks at an arbitrary monad \trip{T}. It is shown that - up to constants - \cat{C} is equivalent to \coalg{K} through the comparison functor (\citep[Theorem 11]{Bar69}). Since \cat{Set} is complete, Theorem \ref{theorem: equivalence between T-algebras and T1-algebras} applies, confirming the general result that the sequence stabilises. However, in this case our result is weaker than the one obtained in \citep{Bar69} by M. Barr. Indeed he showed that \powcat{F}{T} is essentially comonadic as well, whereas Theorem \ref{theorem: equivalence between T-algebras and T1-algebras} does not necessarily imply this, as illustrated by Example \ref{example: ideal functor and downset functor}. At the same time, Example \ref{example: ideal functor and downset functor} indicates that it is not easy, if not impossible, to borrow the assumptions and method used by M. Barr and adapt them for categories enriched over partially ordered sets.\\
\end{example}

\section{Projectives in the category of algebras}
\label{Projectives in the category of algebras}
One of the most interesting aspects of Example \ref{example: ideal functor and downset functor} is the characterisations of projectives in the ambient category as shown in \citep[Proposition 5]{BanNie91}, \citep[Proposition 2.2]{Mad91} and \citep[Theorem 3.2]{Dongsheng97}. These characterisations actually depend on the pair of forgetful-free functors as a parameter that renders the ambient category a category of algebras. However this is not explicitly mentioned. Indeed, J. Madden writes in \citep{Mad91} with regard to the characterisations of stably compact frames as retracts of coherent frames (\cite{Joh81}) that {\em I am sure Johnstone himself has long been aware of this interpretation of his result, but it seems to have escaped many others.} We give here a definition of projectives that varies according to the monad \trip{T} and that specialises into the specific results obtained in the above-mentioned sources. This definition is an easy dualisation of the notion of injective spaces in the sense of Escard\'o (\citep{Esc98,Esc2002}).\\

\begin{definition}
\label{definition: projectives}
Let \trip{T} be a monad on a category \cat{C}. An algebra $(X,a)$ is projective in \alg{T} or simply projective if  $(X,a)$ is projective with respect to the (co)free functor $F^{\mathbb{K}}$.\\
\end{definition}
\noindent
The above definition subsumes Definition 2.1 of \citep{Mad91}, the definition of $E$-projective $Z$-frames in \citep{Dongsheng97} and the definition of projective frames in \citep{BanNie91}. With regard to lax idempotent monads, we have the following equivalences.\\

\begin{theorem}
\label{theorem: dual characterisation of projectives for lax idempotent}
(\citep{Hofmann_Sousa2017,Esc98,Esc2002}) Let \trip{T} be a lax idempotent monad on a category \cat{C}. The following are equivalent:
\begin{enumerate}
\item $(X,a)$ is projective with respect to \set{\map{\alpha}{(TA,m_A)}{(A,\alpha)}\ |\ (A,\alpha)\in\alg{T}}.
\item $(X,a)$ has a structure of a \trip{K}-coalgebra.
\item $(X,a)$ is a retract in \alg{T} of free \trip{T}-algebras. 
\item $(X,a)$ is a \trip{T}-split algebra. (Cf. \citep[Definition 3.6]{Hofmann_Sousa2017})
\end{enumerate}
\end{theorem}
\trip{K}-coalgebras are therefore projectives. In the absence of the condition of lax idempotency, we only have the following characterisation.\\

\begin{theorem}
\label{theorem: dual characterisation of projectives in general}
Projectives are precisely the retracts in \alg{T} of free \trip{T}-algebras.\\
\end{theorem}

\begin{example}
\label{example: projectives for lax idempotent monads}
The projectives in the particular instances of Example \ref{example: ideal functor and downset functor} are respectively: the stably compact frames (\citep{Joh81,Mad91,Ban81}) and the stably continuous $\sigma$-frames (\citep{Mad91,BanMat2003}), the stably supercontinuous frames (\citep[Proposition 5]{BanNie91}), the stably $Z$-continuous semilattices (\citep[Theorem 3.2]{Dongsheng97}).\\
\end{example}

\begin{example}
\label{example: Gleason theorem}
(Gleason, 1958) \citep{Glea58,EilMoo65,Joh82} If \trip{T} is the monad generated by the forgetful functor from the category \cat{KHaus} of compact Hausdorff spaces to \cat{Set} and the free functor \map{\beta}{\cat{Set}}{\cat{KHaus}} - which is the Stone-\v{C}ech compactification of a discrete space, then the projectives with respect to \trip{T} are exactly the extremally disconnected spaces in \cat{KHaus}. Indeed, these are precisely the retracts of spaces of the form $\beta X$. 
\end{example}

\bmhead{Acknowledgements}

I thank B. Jacobs for brief and useful communication with regard to Theorem 4.5 of his paper \citep{Bart2011}.

\section*{Declarations}

\begin{itemize}
\item No funding was obtained for this study.
\item The author declares that there is no conflict of interests.
\item Availability of data and materials not available.
\end{itemize}



\bibliography{sn-bibliography-ideal}

\end{document}